\mag=\magstephalf
\pageno=1
\input amstex
\documentstyle{amsppt}
\TagsOnRight
\interlinepenalty=1000
\NoRunningHeads

\pagewidth{16.5 truecm}
\pageheight{23.0 truecm}
\vcorrection{-1.0cm}
\hcorrection{-1.2cm}
\nologo

\NoBlackBoxes
\font\twobf=cmbx12

\define \dd{\roman d}

\define \sn{{\roman{sn}}}
\define \cn{{\roman{cn}}}
\define \dn{{\roman{dn}}}

\define \ft{{u}}

\define \tvskip{\vskip 1.0 cm}
\define\ce#1{\lceil#1\rceil}
\define\dg#1{(d^{\circ}\geq#1)}
\define\Dg#1#2{(d^{\circ}(#1)\geq#2)}
\define\dint{\dsize\int}
\def\fp{\flushpar}

\define\s#1{\sigma_{#1}}
\define\tp#1{\negthinspace\left.\ ^t#1\right.}
\define\mrm#1{\text{\rm#1}}
\define\lr#1{^{\sssize\left(#1\right)}}

\redefine\qed{\hbox{\vrule height6pt width3pt depth0pt}}

{\centerline{\bf{Closed Loop Solitons and Sigma Functions:   }}}

{\centerline{\bf{Classical and Quantized
 Elasticas with Genera One and Two  }}}

\author
Shigeki MATSUTANI
\endauthor
\affil
8-21-1 Higashi-Linkan Sagamihara 228-0811 Japan
\endaffil \endtopmatter


\centerline{\twobf Abstract }\vskip 0.5 cm

Closed loop solitons in a plane, whose curvatures obey
the modified Korteweg-de Vries equation,
 were investigated. It was shown that
their tangential vectors
 are expressed by ratio of Weierstrass
 sigma functions  for genus one case and ratio of
Baker's sigma functions for the genus two case.
This study is closely related to classical and
quantized elastica problems.


{\centerline{\bf{2000 MSC: 35Q53, 53A04, 14H45, 11J89 }}}

{\centerline{\bf{PACS: 02.30.-f,
05.20.-y,
46.70.Hg,
05.45.Yv
 }}}

 {\centerline{\bf{key words: hyperelliptic function, 
sigma function, space curve,
elastica}}}





\vskip 0.5 cm

\newpage

\centerline{\twobf \S 1. Introduction }\vskip 0.5 cm

This article is on loop solitons of genera one and two.
We give them as  ratios of sigma functions;
for the case of genus one, the sigma functions are
Weierstrass sigma functions while for the case of
genus two, they are Baker's hyperelliptic sigma
functions.

First we mention background of the study, which is
 related to elastica
problem.
Elastica problem was proposed by James Bernoulli in
ending of seventeenth
century
and Euler and his nephew Daniel Bernoulli investigated the
 problem [T, L, W].
The elastica problem is to determine  all allowed shapes
of a thin elastic beam in a plane.
The allowed shape is realized as a curve with the
local minimum of the  energy
$$
        E=\int \dd s k^2, \tag 1-1
$$
where $s$ is its arclength and $k$ is its curvature.
For its tangential angle $\phi$ and $g:= \exp(\sqrt{-1} \phi)$,
the integrand is given by
$$
        k^2 \dd s = -(g^{-1} \dd g) (* g^{-1} \dd g)
        = (\partial_s \phi)^2 \dd s, \tag 1-2
$$
where $*$ is the Hodge star, a map from one-form to zero-form,
and $\partial_s := d /ds$. Thus
the elastica problem is the oldest harmonic map (sigma-model)
problem from a circle $S^1$ or a real line $\Bbb R$ to U(1).
This problem was essentially solved by Euler
using numerical computations [T,L].

For a given elastica in  two dimensional euclidean space,
 it is characterized by
the affine coordinate $(X^1(s), X^2(s))$ or
 $Z(s)=X^1(s)+\sqrt{-1}X^2(s)$. Then its curvature $k$
 is given by $k:=\dfrac{1}{\sqrt{-1}}\partial_s \log \partial_s
 Z$.
The shape locally minimizing the energy (1-1) is given by
the differential equation  [Ma1, Ma2],
$$
        C \partial_s k + \frac{3}{2} k^2 \partial_s k +
         \partial_s^3 k =0.
            \tag 1-3
$$
 This is known as the static modified
Korteweg-de Vries (SMKdV) equation.

As I showed in ref.[Ma2],
when we deal with a non-stretching closed elastica in thermal
bath or investigate the partition function of temperature
 $1/\beta$,
$$
        \Cal Z = \int \Cal D Z \exp( -\beta \int \dd s k^2),
         \tag 1-4
$$
we also need information of curves which have the excited energy.
We will call such curves "quantized
 elastica" following [Ma4]. This problem is associated with
 large polymer problem like DNA in biology and polymer physics.
As shown in [Ma2,Ma3,Ma4], the quantized elastica
 is classified by
 the modified Korteweg-de Vries (MKdV) equation,
$$
        C \partial_t k + \frac{3}{2} k^2 \partial_s k
         + \partial_s^3 k =0,
          \tag 1-5
$$
where
 $\partial_s :=\partial /\partial s$, $\partial_t
 :=\partial /\partial t$ and $t$ is a deformation parameter.
The orbits obeying the MKdV equation preserve local length and
their first integral (1-1), {\it i.e.}, $[\partial_s, \partial_t]=0$
and $\partial_t E=0$.
Since the MKdV equation is an initial value problem,
for a given smooth curve,
there exists, at least, an MKdV orbit including the curve;
 we can classify quantized elastica using the orbits.
Hence  the partition function (1-4) is reduced to
an integration,
$$
        \Cal Z = \int \dd \Omega(E) \exp(-\beta E), \tag 1-6
$$
where $\dd \Omega(E)$ is the density of states.
The evaluation of $\dd \Omega(E)$ means the determination of
the orbits of  the MKdV equation whose first integral is $E$.
In other words, in order to evaluate (1-5),
we must explicitly construct the
periodic solutions of MKdV equation, which are expressed by
hyperelliptic functions.  The hyperelliptic function is
 characterized by genus.
As mentioned in [Ma2], orbits with small genus is expected
to largely contributes to the partition functions (1-4). Thus as
an attempt, we will investigate the classical and quantized
 elastica of
genera one and two respectively.

Since the curve whose curvature obeys the MKdV equation
was called  loop soliton [WIK, I],
 our investigation is, just, a study of  closed loop solitons.
Here we will also note that the time development of
classical elastica does not obey the MKdV equation
except genus one solution case [Ma1].
(From physical point of view, investigation of
a loop soliton with higher genus is somewhat
nonsense in pure kinematic study.)

In this article, we will show that
the tangential vectors of closed loop solitons
 with genera one and two are expressed by ratio of Weierstrass
  sigma functions in \S 3 and
ratio of Baker's sigma functions \S 4 respectively.
This result is very interesting as we will discuss in \S 5.
The sigma function, which is well-tuned theta function,
plays very important roles in the studies of elliptic and
hyperelliptic functions [B1, B2, B3, BEL, G, Ma5, \^O1, \^O2, \^O3].

Further our result has an effect on historical interpretation of
history of science.
According to Truesdell and Love [T,L],
 Euler did illustrate the shape of loop soliton with genus one,
 our result
means that Euler essentially found the elliptic functions and
 sigma functions
even though he dealt only with elliptic integrations.
In fact, the tangential vector $\partial_s Z(s)$ is periodic
 function of $s$.
I believe that Euler noticed its properties and thus after he
read
Fagnano's article, his mind felt the addition formula [W].

The construction of closed loop soliton with
genus two is based upon Baker's hyperelliptic function
theory [B1,B2,B3].
Recently Baker's approach is reevaluated from various
viewpoints [BEL, G,
Ma5, \^O1, \^O2, \^O3]. Buchstaber, Enolskii and Leykin has been
actively studying the
Baker's theory from soliton theory and developing it.
Grant and \^Onishi
found interesting relations in hyperelliptic function in terms of
Baker's theory from point of view of number theory.
As on Baker's theory of the hyperelliptic function,
we start from a hyperelliptic curve and find an explicit function
form without ambiguous parameter, it is very effective.

\tvskip
\centerline{\twobf \S 2. Basic Properties of Loop Soliton }\tvskip

Let us consider a smooth immersion of a circle $S^1$ into the two
dimensional euclidean space $\Bbb E^2 \approx \Bbb C$. The immersed
curve $C$ is characterized by the affine coordinate
$(X^1(s),X^2(s))$.
Here $s$ is a parameter of $S^1$ and is, now,
chosen as the arclength
so that $ds^2 = (dX^1)^2 + (d X^2)^2$.
We will also use the complex expression,
$$
        Z(s):= X^1(s)+ \sqrt{-1} X^2. \tag 2-1
$$
Then  $|\partial_s Z(s)|=1$ and the curvature of $C$ is given by
$$
      k(s) = \frac{1}{\sqrt{-1}} \partial_s \log \partial_s
      Z(s). \tag 2-2
$$

\proclaim{\fp Definition 2-1}\it

A one parameter family of curves $\{C_t\}$ for real parameter
$t\in \Bbb R$
is called a loop soliton,
if its curvature  obeys the MKdV equation;
for $q:=k/2$,
$$
        \partial_t q + 6 q^2 \partial_s q + \partial_s^3 q
         =0. \tag 2-3
$$

\endproclaim

Here we will describe the Miura map for later convenience.

\proclaim{\fp Proposition 2-2 [DJ] }\it

For the solutions $p_\pm$ of the KdV equation,
$$
        \partial_t p_\pm + 6 p_\pm \partial_s p_\pm
        + \partial_s^3 p_\pm =0,
        \tag 2-4
$$
if we find a quantity $q$ satisfying
$p_\pm \equiv q^2 \pm \partial_s q $ over $C_t$ and for all $t$,
$q$ obeys the MKdV equation (2-3).

\endproclaim

\tvskip
\centerline{\twobf \S 3. Loop Soliton with Genus One
and Weierstrass Sigma Function }
\tvskip

In this section, we will deal with a loop soliton with genus one.
In other words, we are consider the case of $t=s/C$,
{\it i.e.}, the SMKdV equation,
$$
        C\partial_s q + 6 q^2 \partial_s q + \partial_s^3 q =0.
        \tag 3-1
$$

First we will set up the Weierstrass $\wp$-functions and
 $\sigma$-functions.

\proclaim{\fp Definition 3-1 [WW]}\it

For an elliptic curve,
$$
       \split
        y^2 &= 4 x^3 - g_2 x + g_3\\
            &= 4(x-e_1)(x-e_2)(x-e_3),
         \endsplit \tag 3-2
$$

\roster

\item We will define  integrals $u$, $\omega_a$, $\eta_a$ ($a=1,2,3$),
$$
        u:= \int^{(x,y)}_{(0,0)} \frac{d x}{y},
        \quad \omega_a = \int^{(e_a,0)}_{(0,0)} \frac{d x}{y},
        \quad \eta_a = \int^{(e_a,0)}_{(0,0)} \frac{x d x}{y}. 
\tag 3-3
$$

\item The elliptic theta function of $\tau := \omega_3/\omega_1$
is defined by
$$
        \theta_1(z) :=  \sum_{n \in \Bbb Z}
        \exp\left( 2\sqrt{-1} \pi
               \left[\frac{1}{2}\tau \left(n-\frac{1}{2}\right)^2
            +(n-\frac{1}{2})( z-\frac{1}{2})+\frac{1}{4}\right]
             \right).
           \tag 3-4
$$

\item
The Weierstrass $\sigma$ function is defined by
$$
        \sigma(u) :=  \exp( \eta_1 u^2/(2\omega_1) )
                   \frac{\theta_1(u/(2\omega_1))}
            {\partial_u \theta_1(u/(2\omega_1))|_{u=0}}. \tag 3-5
$$

\item
We also introduce the other $\sigma$-functions,
$$
        \sigma_a(u):= \exp(-\eta_a u )
              \frac{\sigma(u + \omega_a)}
              {\sigma(\omega_a)}, \quad a=1,2,3. \tag 3-6
$$

\item The Weierstrass $\wp$-function is defined by
$$
          \wp := - \partial_u^2 \log \sigma(u). \tag 3-7
$$

\endroster
\endproclaim

\proclaim{\fp Theorem 3-2 }\it

The shape of loop soliton with genus one  is given by
$$
         Z(s) =\int^s ds \left(\frac{\sigma_3(s-\omega_3/2+\delta)}
       {\sigma(s-\omega_3/2+\delta)}\right)^2, \tag 3-8
$$
where $\delta$ is a constant parameter.
\endproclaim

\tvskip

Later part in  this section devotes the proof of this theorem but
first we will summary the properties of Weierstrass (elliptic)
$\wp$ and $\sigma$ functions.
Let "$'$" be derivative in $u$.

\proclaim{\fp Proposition 3-3 [WW]} \it

\roster

\item
$$
        \split
        \wp'(u)^2&= 4 \wp^3 + g_3 \wp +g_4 \\
                           &=4(\wp-e_1) (\wp -e_2) (\wp- e_3),
         \endsplit \tag 3-9
$$
where $\wp(\omega_a) = e_a$ and $e_1+e_2+e_3=0$.

\item
$$
        12 \wp(u) \partial_u \wp(u)+ \partial_u^3 \wp(u)=0.
         \tag 3-10
$$

\endroster
\endproclaim

Further we will note the addition relations.

\proclaim{\fp Proposition 3-4 [WW]} \it

\roster

\item
$$
        \wp(z) - \wp(u) = \frac{\sigma(z+u)\sigma(z-u)}
        {[\sigma(z)\sigma(u)]^2}.          \tag 3-11
$$

\item
$$
        \wp(u+z) - \wp(u)
 = -\frac{1}{2}\partial_u\frac{\wp'(z)-\wp'(u)}{\wp(z)-\wp(u)}.
  \tag 3-12
$$

\item
$$
        \wp(z+u) + \wp(z) +\wp(u)
           = \frac{1}{4}\left(\frac{\wp'(u)-\wp'(z)}
           {\wp(u)-\wp(z)}\right)^2.
          \tag 3-13
$$

\item
$$
\wp(u)-e_a = \left(\frac{\sigma_a(u)}{\sigma(u)}\right)^2.
 \tag 3-14
$$
\endroster
\endproclaim

Here (3-12) and (3-14) is directly derived from (3-11).
(3-13) is slightly difficult but using proposition 3-3, it can
be proved.

\tvskip

From the definition of loop soliton, the theorem 3-2 is
reduced to the next lemma.

\proclaim{\fp Lemma 3-5}\it

For $\rho:=\wp(u) -e_3$, $\mu(u):=\dfrac{1}{2\sqrt{-1}}
 \partial_u \log \rho$
obeys the MKdV equation,
$$
        \partial_u \mu + 6 \mu^2 \partial_u \mu
         + \partial_u^3 \mu =0.
       \tag 3-15
$$
\endproclaim

\demo{Proof}
From proposition 3-4,
$$
        \wp(u+\omega_3) - \wp(u) =
        \frac{1}{\sqrt{-1}}\partial_u\mu,\tag 3-16
$$
$$
        \wp(u+\omega_3) + \wp(u) +e_3
           = -\mu^2. \tag 3-17
$$
Thus for $v:=u -\omega_3/2$, we have a relation,
$$
        -2 \wp(v \pm \omega_3/2)- e_3
           = \mu^2(v)\pm \sqrt{-1}\partial_v \mu(v).
$$
From (3-10),  $p=-2\wp(u)-e_3$ satisfies the equation,
$$
        6 e_3 \partial_u p+6 p \partial_u p + \partial_u^3 p =0.
           \tag 3-18
$$
Using the proposition 2-2, this lemma is proved. \qed \enddemo

Hence the theorem 3-2 is also satisfied.

\proclaim{\fp Remark 3-6}

We note that (3-8) is a necessary condition.
If one tunes $\delta$ in (3-8) and the $g_2$ and $g_3$
in (3-2) to find a real line $s$ in the complex
plane $\{u\}= \Bbb C$ so that it satisfies the reality condition
$|\partial_s Z|=1$ and the closed condition $Z(s)=Z(s+2\omega_1)$,
(3-8) becomes the closed loop soliton.
\endproclaim

\tvskip
\centerline{\twobf \S 4. Loop Soliton with Genus Two
 and Baker's Sigma Function }\tvskip

In  this section, we will deal with  a hyperelliptic solution
of the loop soliton whose genus is two.
We will redefine the quantities
$\sigma$, $\wp$, $\rho$,
$\omega$, $\eta$, $\tau$ in the previous section.
These quantities are altered to  genus two versions of corresponding
quantities in the previous section.

\proclaim {\fp Definition 4-1 [B1, B2, BEL
Chapter 2, \^O1 p.384-386, \^O2]}\it

For a hyperelliptic curve   of genus two,
$$
 \split
   y^2 &= f(x) \\
   &= \lambda_0 +\lambda_1 x
        +\lambda_2 x^2  +\cdots +\lambda_5 x^{5},
\endsplit \tag 4-1
 $$
where $\lambda_{5}\equiv1$ and $\lambda_j$'s are complex numbers,
we will give definitions as follows.

 \roster

\item
We choose $a_r$, $c_r$ and $c$ $(r=1,2)$ for
$$
   f(x) =P(x)Q(x),
$$ $$
 P(x)=(x-a_1)(x-a_2),\quad Q(x)=(x-c_1)
     (x-c_2)(x-c), \tag 4-2
$$
so that curves $(a_r,c_r)$ or $(c,\infty)$ does not
intersect each other and
 $Re (a_r)\le Re (c_r)$, $Re (a_1) \le Re (a_{2})$,
$Re (c_2) \le Re (c)$.

\item  Let us denote the homology of the hyperelliptic
curve  by
$$
\roman{H}_1(X_g, \Bbb Z)
  =\Bbb Z\alpha_{1}\oplus\Bbb Z\beta_{1}
    \oplus\Bbb Z\alpha_{2}\oplus\Bbb Z\beta_{2},\tag 4-3
$$
where these intersections are given as
$[\alpha_i, \alpha_j]=0$, $[\beta_i, \beta_j]=0$ and
$[\alpha_i, \beta_j]=\delta_{i,j}$.

\item The unnormalized  differentials of first kind are defined by,
$$   d u_1 := \frac{ d x}{2y}, \quad
      d u_2 :=  \frac{x d x}{2y}.\tag 4-4
$$

\item The unnormalized differentials of second kind are defined by,
$$
     d r_{1}:=\dfrac{1}{2y}\left(
      \lambda_3 x+ 2\lambda_4 x^2+ 3\lambda_5 x^3 \right)dx ,\quad
     d r_{2}:=\dfrac{1}{2y}
      \lambda_{5} x^2 dx .\tag 4-5
$$

\item The unnormalized period matrices are defined by,
$$   2 \pmb{\omega}':=\left[\matrix
        \dint_{\alpha_{1}}d u_{1} &  \dint_{\alpha_{2}}d u_{1} \\
         \dint_{\alpha_{1}}d u_{2} &  \dint_{\alpha_{2}}d u_{2}
         \endmatrix \right],
      2\pmb{\omega}'':=\left[\matrix
        \dint_{\beta_{1}}d u_{1} &  \dint_{\beta_{2}}d u_{1} \\
         \dint_{\beta_{1}}d u_{2} &  \dint_{\beta_{2}}d u_{2}
         \endmatrix \right],
 \quad
    \pmb{\omega}:=\left[\matrix \pmb{\omega}' \\ \pmb{\omega}''
     \endmatrix\right].\tag 4-6
$$

\item The normalized period matrices are given by,
$$
   \pmb \tau:={\pmb{\omega}'}^{-1}\pmb{\omega}'',
   \quad
    \hat{\pmb{\omega}}:=\left[\matrix 1_g \\ \pmb \tau
     \endmatrix\right].\tag 4-7
$$

\item The complete hyperelliptic integral of the second kinds
is given  as
$$
   2 \pmb{\eta}':=\left[\matrix
        \dint_{\alpha_{1}}d r_{1} &  \dint_{\alpha_{2}}d r_{1} \\
         \dint_{\alpha_{1}}d r_{2} &  \dint_{\alpha_{2}}d r_{2}
         \endmatrix \right],
      2\pmb{\eta}'':=\left[\matrix
        \dint_{\beta_{1}}d r_{1} &  \dint_{\beta_{2}}d r_{1} \\
         \dint_{\beta_{1}}d r_{2} &  \dint_{\beta_{2}}d r_{2}
         \endmatrix \right].
      \tag 4-8
$$

\item We will defined the Riemann theta function over $\Bbb C^2$
characterized by $\hat{ \pmb{\Lambda}}$,
$$
   \split
   \theta\negthinspace\left[\matrix a \\ b \endmatrix\right](z)
   & :=\theta\negthinspace\left[\matrix a \\ b \endmatrix\right]
     (z; \pmb\tau) \\
   & :=\sum_{n \in \Bbb Z^2} \exp \left[2\pi \sqrt{-1}\left\{
    \dfrac 12 \ ^t\negthinspace (n+a)\pmb \tau(n+a)
    + \ ^t\negthinspace (n+a)(z+b)\right\}\right],
 \endsplit
         \tag 4-9
$$
for $2$-dimensional vectors $a$ and $b$.

\item We will introduce constant parameters, [B2 p.343]
$$
      \lambda_r = \frac{(-1)^{r} \sqrt{P'(a_r)}}
                       {(\sqrt{-1f'(a_r)/4})^{1/2}}, \quad
        r=1,2.
    \tag 4-10
$$
\endroster
\endproclaim

\tvskip
We will note that these contours in the integral are,
for example,
given in p.3.83 in [M]. Thus above values can be, in principle,
computed in terms of numerical method for a given $y^2=f(x)$.

\vskip 0.5 cm
\proclaim {Proposition 4-2 [B1, B2 p.316-p.335,
BEL p.25, M p.3.80-84, \^O2]}

The Riemannian constant $K\in \Bbb C^2$,
$\omega_a \in \Bbb C^2$ are given by
$$
 K = \pmb \omega^{\prime -1}(\int^{(a_1,0)}_\infty
  d\bold u+\int^{(a_2,0)}_\infty d\bold u )
=  \delta' + \delta'' \pmb\tau, \tag 4-11
$$
$$
\split
\omega_1:=\int^{(a_1,0)}_\infty d\bold u&
      = \pmb \omega^{\prime}\delta_1'
      +\pmb\omega^{\prime\prime}\delta_1'',\\
\omega_2:=\int^{(a_2,0)}_\infty d\bold u&
       = \pmb \omega^{\prime}\delta_2'
       +\pmb\omega^{\prime\prime}\delta_2'',
\endsplit
\tag 4-12
$$
where
$$
 \delta'_1 :=\left[\matrix \dfrac {1}{2}\\
          {}\\ 0\endmatrix\right],\quad
 \delta'_2 :=\left[\matrix \dfrac {1}{2}\\
       \dfrac {1}{2}\endmatrix\right],\quad
  \delta' =\left[\matrix 1 \\ {}\\
       \dfrac {1}{2}\endmatrix\right],
 $$
$$
 \delta''_1:=\left[\matrix \dfrac{1}{2} \\
           {}\\ 0\endmatrix\right],\quad
 \delta''_2:=\left[\matrix 0\\ {}\\
          \dfrac{1}{2}  \endmatrix\right],\quad
 \delta'':=\left[\matrix \dfrac{1}{2} \\
          \dfrac{1}{2} \endmatrix\right].
\tag 4-13
$$
\endproclaim

\tvskip

\vskip 0.5 cm
\proclaim {Definition 4-3 [B1 p.286, B2 p.336, p.353, p.370,
 BEL p.32, p.35, \^O1 p.386-7, \^O2] }

\it
\roster
We will define the coordinate $(u_1, u_2)$ in $\Bbb C^2$
for  points $(x_1,y_1)$ and $(x_2,y_2)$
of the curve $y^2 = f(x)$,
$$
  \ft_j :=\int^{(x_1,y_1)}_\infty du_j
       +\int^{(x_2,y_2)}_\infty du_j .\tag 4-14
$$

\item Using the coordinate $\ft_j$, sigma functions,
which are homomorphic
functions over $\Bbb C^2$, are defined by,
$$ \sigma(\ft)=\sigma(\ft;\pmb\omega):
  =\gamma \roman{exp}(-\dfrac{1}{2}\ ^t\ \ft
  \pmb\eta'{\pmb{\omega}'}^{-1}\ft)
  \vartheta\negthinspace
  \left[\matrix \delta'' \\ \delta' \endmatrix\right]
  ({(2\pmb{\omega}')}^{-1}\ft ;\Bbb T) , \tag 4-15
$$
$$
 \split
 \sigma_r(\ft)&=\sigma_r(\ft;\pmb\omega):
  =\lambda_r  \roman{exp}(-\dfrac{1}{2}\ ^t\ \ft
  \pmb \eta'{\pmb{\omega}'}^{-1}\ft)
  \vartheta\negthinspace
  \left[\matrix \delta_r'' \\ \delta_r' \endmatrix\right]
  ({(2\pmb{\omega}')}^{-1}\ft ;\Bbb T),\\
  & \qquad \qquad \text{for }a=0,1,2,
  \endsplit \tag 4-16
$$
where $\gamma$ is a fixed constant,
$\delta_0'\equiv \delta''_0:=\left[\matrix 0 \\
            0\endmatrix\right]$, $\lambda_0=1$,

\item
In terms of $\sigma$ function, $\wp$-functions over the
hyperelliptic curve are given by
$$   \wp_{\mu\nu}(\ft)=-\dfrac{\partial^2}{\partial
   \ft_\mu\partial \ft_\nu}
   \log \sigma(u) ,
$$
$$   \wp_{\mu\nu}^{(r)}(\ft)=-\dfrac{\partial^2}{\partial
   \ft_\mu\partial \ft_\nu}
   \log \sigma_r(u) . \tag 4-17
$$
\endroster

\endproclaim

\proclaim {\fp Theorem 4-4} \it

The shape of loop soliton with genus two is given by
$$
         Z(s,t) =\int^s ds
           \left( \frac{\sigma_2(u(s,t)+\omega_1+\omega_2/2+\delta)}
           {\sigma_0(u(s,t)+\omega_1+\omega_2/2+\delta)}\right)^2,
       \tag 4-18
$$
for $a_2 =-\dfrac{1}{3} \lambda_4$, where $\delta$ is
a constant parameter and $u(s,t):=\pmatrix -4 t \\ s \endpmatrix$.

\endproclaim

\tvskip

Before proving this theorem, we will also review
the properties of Baker's $\wp$ and $\sigma$-functions.

\proclaim{\fp Proposition 4-5 [B3 p.155-6, BEL p.56-58,
 G p.103-4, p.116, \^O1 p.388,\^O2}\it

Let us express
$\wp_{\mu\nu\rho}:=\partial \wp_{\mu\nu}(\ft)
  /\partial \ft_\rho$ and
$\wp_{\mu\nu\rho\lambda}:=\partial^2
 \wp_{\mu\nu}(\ft) /\partial \ft_\mu \partial \ft_\nu$.
Then hyperelliptic $\wp$-functions obey the relations
$$   \allowdisplaybreaks \align
     (H-1)\quad& \wp_{2222}-6\wp_{22}^2
      = 2 \lambda_3 \lambda_5
           + 4 \lambda_4 \wp_{22}
       + 4\lambda_5 \wp_{21}
         ,\\
     (H-2)\quad& \wp_{2221}-6\wp_{22} \wp_{21}
      =
           4 \lambda_4 \wp_{21}
       - 2\lambda_5 \wp_{11}, \\
     (H-3)\quad& \wp_{2211}-4\wp_{21}^2
            -2\wp_{22} \wp_{11}
      = 2\lambda_3 \wp_{21}, \\
     (H-4)\quad& \wp_{2111}-6\wp_{21} \wp_{11}
      =- 4\lambda_0 \lambda_5
      - 2 \lambda_1\wp_{22}
       + 4\lambda_2 \wp_{21}, \\
     (H-5)\quad& \wp_{1111}-6\wp_{11}^2
      =- 8\lambda_0 \lambda_4+ 2 \lambda_1 \lambda_3
       -12\lambda_0 \wp_{22}
       + 4\lambda_1 \wp_{21} + 4 \lambda_2 \wp_{11}.
     \endalign
$$
$$   \allowdisplaybreaks \align
     (I-1)\quad& \wp_{222}^2
      = 4(\wp_{22}^3 +\wp_{12}\wp_{22}+\lambda_4\wp_{22}^2+
            \wp_{11}+\lambda_3\wp_{22} +\lambda_2)
         ,\\
     (I-2)\quad& \wp_{222}\wp_{221}
      =4 (\wp_{12}\wp_{22}^2 -\frac{1}{2}(\wp_{11}\wp_{22}
             -\wp_{12}^2\\
             &\qquad + \lambda_3 \wp_{12}-\lambda_1) +
              \lambda_3\wp_{12} +\lambda_4 \wp_{12}\wp_{22}),\\
     (I-3)\quad& \wp_{221}^2
      =
          4\left (\wp_{11}\wp_{22}^2 -
              (\wp_{11}\wp_{22} - \wp_{12}^2
                 + \lambda_3 \wp_{12}-\lambda_1)\wp_{22}
             - \wp_{11}\wp_{12} \right.\\
            & \qquad + \lambda_4 \wp_{11}\wp_{22}
                  + \lambda_3 \wp_{12} \wp_{22}
             -\lambda_4 (\wp_{11}\wp_{22}-\wp_{12}^2
              +\lambda_3 \wp_{12}
             -\lambda_1) \\
          & \qquad\left. + \lambda_4 \lambda_3 \wp_{12}
             -\lambda_1 \wp_{22}
              -\lambda_1\lambda_4+\lambda_0 \right).
     \endalign
$$

\endproclaim

We also have the additional formulae for them

\proclaim{\fp Proposition 4-6 [B2 p.372, p.381, BEL p.112]}\it

\roster

\item
$$
\frac{\sigma(u+v)\sigma(u-v)}{\sigma(u)^2\sigma(v)^2}
 =\wp_{22}(u)\wp_{21}(v)-\wp_{21}(u)\wp_{22}(v) -\wp_{11}(u)
 +\wp_{11}(v). \tag 4-19
$$

\item
$$
        \left(\frac{\sigma_2(u)}{\sigma_0(u)}\right)^2
          = \frac{1}{Q(a_r)}
           (\wp^{(0)}_{21}(u)+\wp^{(0)}_{22}(u)a_2-a_2^2).
           \tag 4-20
$$
\endroster

\endproclaim

We note that (4-19) corresponds to (3-11) and (4-20)
is a genus two version of
(3-14). The similar formulae of (3-13) and (3-14)
are studied in [G p.114-6]
but they are not so effective for this purpose. For special value 
(theta
characteristics), we give genus two versions (4-24) and (4-26)
of (3-13) and (3-14) in proof of lemma 3-5.

\proclaim{\fp Proposition 4-7}\it

$$
        \wp^{(0)}_{\mu\nu}(u)=\wp_{\mu\nu}(u-\omega_1-\omega_2),
$$
$$
        \wp^{(1)}_{\mu\nu}(u)=\wp_{\mu\nu}(u-\omega_2), \quad
        \wp^{(2)}_{\mu\nu}(u)=\wp_{\mu\nu}(u-\omega_1), \quad
           \tag 4-21
$$

\endproclaim
\demo{Proof}
Using $\theta(z):=
  \vartheta\negthinspace
  \left[\matrix \delta_0'' \\ \delta_0' \endmatrix\right](z)$,
the $\sigma$-functions are explicitly written by
$$
      \split
        \sigma( u ) &=\gamma
        \exp\left(
        -\frac{1}{2} {}^t u \pmb\eta'\pmb\omega^{\prime -1}
         u +2 \pi \sqrt{-1}
         \left[ \frac{1}{2}{}^t\delta'' \pmb\tau \delta''
+ ^t\delta''({(2\pmb\omega')}^{-1}u + \delta')\right] \right) \\
      &\qquad
    \times \vartheta({(2\pmb\omega'}^{-1}u + \delta'+\pmb\tau \delta'') ,
    \endsplit
$$
$$
     \split
        \sigma_r( u ) &=\lambda_r
        \exp\left(
        -\frac{1}{2} {}^t u \pmb\eta'\pmb\omega^{\prime -1} u
         +2 \pi \sqrt{-1}
         \left[ \frac{1}{2}{}^t\delta_r'' \pmb\tau \delta_r''
+ ^t\delta_r''({(2\pmb\omega')}^{-1}u + \delta_r')\right] \right)  \\
    &\qquad \times
    \vartheta({(2\pmb\omega'}^{-1}u + \delta_r'+\pmb\tau \delta_r'').
    \endsplit \tag 4-22
$$
Thus
$$
 \log(\sigma_r(u)/\sigma(u)) -
  \log( \vartheta({(2\pmb\omega'}^{-1}u + \delta_r'+\pmb\tau \delta_r'')
   / \vartheta({(2\pmb\omega'}^{-1}u + \delta'+\pmb\tau \delta''))
$$
is a first order polynomial of $u$. Thus the derivative in
the definitions of $\wp$ leave only the difference of phases.
Noting proposition 4-2 and $K=(2\pmb\omega')^{-1}(\omega_1+\omega_2)$,
we obtain the identities,
$$
       u=u-\omega_1-\omega_2 + 2\pmb\omega' K,
$$ $$
       u+2\pmb\omega'(\delta_1' +\pmb\tau \delta_1'')
    =u-\omega_2 + 2\pmb\omega' K,
$$ $$
    u+2\pmb\omega'(\delta_2' +\pmb\tau \delta_2'')
    =u-\omega_1 + 2\pmb\omega' K. \tag 4-23
$$
Hence (4-21) is proved.
\qed \enddemo

From the definition 2-1, the theorem 4-4 is reduced to next lemma.

\proclaim{\fp Lemma 4-8}\it

For $a_2 =-\dfrac{1}{3} \lambda_4$,
$\mu(u):=\dfrac{1}{\sqrt{-1}}\partial_{u_2}\log
\left(\dfrac{\sigma_2(u)}{\sigma_0(u)}\right)$
obeys the MKdV equation,
$$
        -4\partial_{u_1} \mu +
         6 \mu^2 \partial_{u_2} \mu +
         \partial_{u_2}^3 \mu =0.
       \tag 4-24
$$
\endproclaim

\demo{Proof}
By noting (4-20) and introducing
$$
        \rho(u):=  \wp^{(0)}_{21}(u)+\wp^{(0)}_{22}(u)a_2-a_2^2,
        \tag 4-25
$$
$\mu(s)$ is given by
$$
        \mu(u) = \frac{1}{2\sqrt{-1}}
            \frac{\partial_{u_2}\rho(u)}{\rho(u)}. \tag 4-26
$$

By taking the logarithm and derivative in $u_2$,
the addition formula (4-20) is reduced to
$$
 \wp_{22}(u+\omega_2) - \wp_{22}(u)
 = \sqrt{-1}\partial_{u_2} \mu(u+\omega_1+\omega_2). \tag 4-27
$$
Thus we must prove the relation
$$
 \wp_{22}(u+\omega_2) + \wp_{22}(u) + \lambda_4 +a_2
          = - \mu(u+\omega_1+\omega_2)^2. \tag 4-28
$$
If (4-28) accomplishes, (4-27) and (4-28) give the
relations
that  for $v:=u-\omega_2/2$,
$$
-2\wp_{22}(v\pm\omega_2/2)-\lambda_4-a_2
         = \mu^2(v+\omega_1+\omega_2) \pm \partial_{u_2}
         \mu(v+\omega_1+\omega_2). \tag 4-29
$$
Using (H-1) in the proposition 4-5,
we can show that $q:=-2\wp_{22}(v\pm\omega_2/2)-\lambda_4-a_2 $ obeys
$$
        -4 \partial_{u_1} q +
         6 q \partial_{u_2} q +
         \partial_{u_2}^3 q =0,
       \tag 4-30
$$
if $a_2= -\lambda_4/3$. From the proposition 2-2,
if (4-28) is true, lemma 4-8 will be proved.

Thus we will concentrate (4-28).
We add (4-27) to $2 \wp(u)+\lambda_4+a_2$ and then
we have the relation
$$
 \wp_{22}(u+\omega_2) + \wp_{22}(u) + \lambda_4 +a_2
          = - \mu(u+\omega_1+\omega_2)^2
          + \frac{\Delta(u+\omega_1+\omega_2)}
          {4 \rho(u+\omega_1+\omega_2)^2},
           \tag 4-31
$$
where,
$$
        \Delta(u) := 8 \wp_{22}(u)\rho(u)^2 - 2(\partial_{u_2}^2
                      \rho(u)) \rho(u)
              +(4 \lambda_4+a_2+1)(\partial_{u_2}\rho_2(u))^2.
               \tag 4-32
$$
We use the relations in proposition 4-5 and then
direct computation shows that $\Delta$ vanishes
even though it is very tedious.
Accordingly the lemma is proved \qed \enddemo

\tvskip

The theorem 4-4 is also proved and the remark 3-6
 is also effective in this case.

\tvskip
\centerline{\twobf \S 5. Discussion }\tvskip

We showed that the loop soliton is directly
connected with $\sigma$-functions.
In ref.[Ma1], I partially showed the fact but
I could not reach the concrete function form
and the difference between numerator and
denominator in the ratio of $\sigma$ or ($\tau$)
function.
This study determines that the difference comes from
the half of period or so-called theta
characteristics [B2, M].
Even though we could not find a solution with
higher genus, it is not difficult to
conjecture the function form of loop soliton
in terms of the theta characteristics.

As we used the definition of loop soliton as
a relation of MKdV equation,  our study
is just that of MKdV equation.
In other words,
we have concretely constructed the function forms
of the genus one and
the genus two periodic solutions of MKdV equation in
the section 3 and 4.
In ordinary study of elliptic function solution
of MKdV equation [Ma1], one uses Jacobi elliptic functions,
$\sn$, $\cn$, and $\dn$ functions.
Our solution (3-15)
becomes the Jacobi elliptic function by means of
Landen transformation, Jacobi transformation,
Gauss transformation and so on.
The relation between
Weierstrass and Jacobi elliptic functions is
derived from (3-14) [WW]. $\sqrt{ \wp(u)-e_a}$ plays
very important roles;
$$
        \sn(z) = \sqrt{\frac{e_1-e_3}{\wp(u)-e_3}},
        \quad
        \cn(z) = \sqrt{\frac{\wp(u)-e_1}{\wp(u)-e_3}},
        \quad
        \dn(z) = \sqrt{\frac{\wp(u)-e_2}{\wp(u)-e_3}},
        \tag 5-1
$$
where $z=(e_1- e_3)^{1/2} u $. Since $\sn$ function
is connected with $1/\sn$ by phase shift and associated
with $\dn$ and $\cn$ functions by algebraic equations,
the Jaobi elliptic function
could roughly be considered  as $\sqrt{ \wp(u)-e_3}$.
Similarly we might regard $\sqrt{\rho}$ in (4-25) as
a Jacobi-type hyperelliptic function due to (4-20);
$\sqrt{\rho/a_2}$ is essentially   $\sqrt{ \wp(u)-e_3}$ when
we set $\partial_{u_1} \sigma = 0$ and $u:=u_2$.

We should also notice  the square in
(3-14) and (4-20). As described in [Ma1 and reference therein],
it is closely related to fermion. In fact,
$\psi:=\sigma_2/\sigma$ $\sim \sqrt{\partial_s Z}$ in (3-14)
( $\psi:=\sigma_2/\sigma_0$ in (4-20) )
obeys the Dirac equation,
$$
        \pmatrix \partial_s & \mu \\
                 \mu & -\partial_s \endpmatrix
                 \pmatrix \psi \\ \sqrt{-1} \psi \endpmatrix =0.
                 \tag 5-2
$$
The square root of $\partial_s Z$ means double covering of
the tangent bundle. (5-2) is related to Schwarz derivative and
moduli space as mentioned in [Ma1].
We recognize that there is still open problem what is the orign of
the square root.

Finally we will comment upon the closed condition
mentioned in the remark 3-6.
We could not touch the condition. In general,
an integration of periodic function is not periodic.
However if the function satisfies some conditions,
its integral is also periodic.
The closed condition is a such condition.
It suppresses two degree of freedom in
the parameter spaces.
For the study of quantized elastica problem,
it is very important to determine the relations
between the conditions and the coefficients
of the algebraic curve $\lambda_j$ in (4-1).
In future, we wish to determine them.

\tvskip
{\centerline{\bf{ Acknowledgment}}}
\tvskip

I thank Prof.~Y.~\^Onishi for leading me this beautiful theory
of Baker. I am grateful to Prof.~K. Tamano and
H.~Mitsuhashi for fruitful discussions.
I am grateful to Prof. V. Z. Enolskii for sending me
his papers on the Abel functions related to
Baker's theory.

\enddocument

\ref
  \key   {\bf {GGKM}}
  \by    Gardner, C.S., Greene J. M., Kruskal M. D.,
         and Miura, R. M.
  \paper Method for Solving the Korteweg-de Vries Equation
  \vol   19
  \yr    1967
  \pages 1095-1097
  \jour  Phys. Rev. Lett.
\endref

\endRefs


\Refs
\widestnumber\key{BBEIM}

\ref
  \key   {\bf {B1}}
  \by    Baker, H.F.
  \book  Abelian functions
         -- Abel's theorem and the allied theory
            including the theory of the theta functions --
  \publ  Cambridge Univ. Press
  \yr    1897, republication 1995
\endref
\ref
  \key   {\bf {B2}}
  \by    Baker, H.F.
  \paper On the hyperelliptic sigma functions
  \jour  Amer. J. of Math.
  \vol   XX
  \yr    1898
  \pages 301-384
\endref
\ref
  \key   {\bf {B3}}
  \by    Baker, H.F.
  \paper On a system of differential equations
leading to periodic functions
  \jour  Acta math.
  \vol   27
  \yr    1903
  \pages 135-156
\endref
\ref
  \key   {\bf {BEL}}
  \by    Buchstaber, V.H., Enolskii, V.Z. and Leykin, D.V.
  \paper Klein Function, Hyperelliptic Jacobians and
         Applications
  \jour Rev. Math. \& Math. Phys.
  \yr    1997 \vol 10 \pages 3-120
\endref

\ref
  \key   {\bf {DJ}}
  \by    Drazin P. G. and Johnson R. S.
  \book  Solitons: an introduction
  \publ  Cambridge Univ. Press
  \yr    1989
\endref

\ref
  \key   {\bf {G}}
  \by    Grant, D.
  \paper Formal groups in genus two
  \vol   411
  \yr    1990
  \pages 96--121
  \jour  J. reine angew. Math.
\endref



\ref
  \key   {\bf {I}}
  \by    Ishimori, Y
  \paper On the Modified Korteweg-de Vries Soliton
    and the Loop Soliton \jour J. Phys. Soc. Jpn
  \yr    1981\vol 50 \pages 2741-2472
\endref


\ref \key {\bf {KIW}}
  \by    Konno, K. Ichikawa, Y. and Wadati M.
  \paper A  Loop Soliton propagating along a Stretched Rope
  \jour J. Phys. Soc. Jpn
  \yr    1981\vol 50
       \pages 1025-1026 \endref
\ref \key {\bf{L}}
 \by A.~E.~H.~Love \book A Treatise on the Mathematical Theory
   of Elasticity \publ Cambridge Univ. Press \yr 1927
   \publaddr Cambridge \endref



\ref
  \key   {\bf {M}}
  \by    Mumford, D.
  \book  Tata lectures on theta II {\rm (Prog. in Math. vol.43)}
  \yr    1984
  \publ  Birkh\"auser
\endref
\ref
  \key   {\bf {Ma1}}
  \by    Matsutani, S.
      \paper Geometrical Construction of the Hirota Bilinear
        Form of the Modified Korteweg-de Vries Equation
        on a Thin Elastic Rod: Bosonic Classical Theory
        \jour Int. J. Mod. Phys. A \yr 1995
          \vol 22 \pages 3109-3123
        \endref
\ref
  \key   {\bf {Ma2}}
      \bysame
      \paper Statistical Mechanics of Elastica on plane:
      Origin of MKdV hierarchy
        \jour J. Phys. A \yr 1998  \vol 31 \pages 2705-2725
        \endref
\ref
  \key   {\bf {Ma3}}
  \bysame
      \paper Statistical Mechanics of Elastica in $\Bbb R^3$
      \jour J. Geom. Phys. \vol 29 \yr 1999 \pages 243-259
        \endref

\ref\key   {\bf {Ma4}}
  \bysame
      \paper On the Moduli of a Quantized Elastica in $\Bbb P$ and
             KdV Flows:
         Study of Hyperelliptic Curves as an Extension of Euler's
         Perspective of Elastica I
      \jour mathDG/9808099
        \endref
\ref\key   {\bf {Ma5}}
  \bysame
      \paper Hyperelliptic Solutions of KdV and
KP equations:Reevaluation of Baker's Study on
Hyperelliptic Sigma Functions
      \jour nlin.SI/0007001
        \endref

\ref \key {\bf{\^O1}} \by \^Onishi Y. \paper Complex
multiplication formulae for curves of genus three
\jour Tokyo J. Math. \vol 21 \pages 381-431 \yr1998
\endref
\ref \key {\bf{\^O2}} \bysame \paper chodaenkansuu-ron
(Introduction to Hyperelliptic Function)
\jour unpublished \lang japanese  \yr1998
\endref
\ref \key {\bf{\^O3}}\bysame \paper Determinatal Expressions
for Some Abelian Functions in Genus Two
\jour preprint \yr2000
\endref

\ref \key {\bf{T}} \by Truesdell, C.
\jour Bull. Amer. Math. Soc. \vol 9 \yr
    1983 \page 293-310 \paper
    The influence of elasticity on analysis:
    the classic heritage \endref
\ref \key {\bf{W}}
 \by Weil, A. \book Number Theory: an approach through history;
From Haammurapi to Legendre  \publ
Birkh\"auser \publaddr Cambridge \yr 1983
\endref

\ref \key {\bf{WW}} \by E.~T.~Wittaker and G.~N.~Watson
\book A Course of Modern Analysis
\publ Cambridge Univ. Press \yr 1927 \publaddr Cambridge \endref


\endRefs
\enddocument